\documentclass[a4paper,twoside,12pt]{article}

\usepackage[T1]{fontenc}

\usepackage{fancyhdr}

\usepackage{a4}
\usepackage{amsmath}
\usepackage{amssymb}
\usepackage{amsthm}
\usepackage{verbatim}
\usepackage[dvips]{graphicx}
\usepackage[utf8x]{inputenc} 
\usepackage[polish,british]{babel} 

\newenvironment{prf}{\par\noindent{\bf Proof.}}{\par\rightline{$\Box$}}

\numberwithin{equation}{section}
\newtheorem{theo}[equation]{Theorem}
\newtheorem{lem}[equation]{Lemma}
\newtheorem{cor}[equation]{Corollary}

\newtheorem{defin}[equation]{Definition}

\newtheorem{ex}[equation]{Example}
\newtheorem{rem}[equation]{Remark}
\newtheorem{prop}[equation]{Proposition}

\newtheorem{notation}[equation]{Notation}

\def\s{\sigma}
\def\t{\tau}
\def\r{\rho}
\def\D{\Delta}

\def\conv{\operatorname{conv}}

\def\cone{\operatorname{cone}}

\def\e{\varepsilon}
\def\s{\sigma}
\def\mono{\hookrightarrow}

\def\C{{\mathbb C}}

\def\N{{\mathbb N}}

\def\R{{\mathbb R}}
\def\Z{{\mathbb Z}}
\def\P{{\mathbb P}}

\def\ra{\rightarrow}

\def\la{\leftarrow}

\def\izo{\simeq}

\def\ccO{{\cal O}}

\def\Spec{\operatorname{Spec}}
\def\Hom{\operatorname{Hom}}
\def\Pic{\operatorname{Pic}}
\def\Sing{\operatorname{Sing}}
\def\codim{\operatorname{codim}}
\def\rk{\operatorname{rk}}

\newcommand{\mr}{\longrightarrow}
\newcommand{\mapdown}[1]{\Big \downarrow
  \rlap{$\vcenter{\hbox{$\scriptstyle#1$}}$}}

\def\Li#1{({#1}_1,\ldots,{#1}_n)}
\def\Lo#1{({#1}_0,\ldots,{#1}_n)}

\def\e{\varepsilon}
\def\p{\varphi}
\def\0{\{0\}}

\author{Weronika Buczy\'n{}ska} 
\title{Fake weighted projective spaces}
\date{June 2002}

\begin{document}

\maketitle

\begin{abstract}
We define fake weighted projective spaces as a generalisation of weighted projective spaces. 
We introduce the notions of fundamental group in codimension 1 
and of universal covering in codimension 1. 
We prove that for every fake weighted projective space its universal cover in codimension 1 
is a weighted projective space.
\end{abstract}

\tableofcontents

\bigskip

\noindent \textbf{Keywords:} toric geometry, fundamental group in codimension 1, 
cover in codimension 1, weighted projective space, fake weighted projective space.

\section{Introduction.}
Miles Reid in his very well known paper \cite{reid_decomp}, 
states (the corollary 2.5 about contraction of extremal rays):

\begin{theo}\label{Reid}
If $\D$ is a complete simplicial fan, $R$ an extremal  ray of $NE(X(\D))$ then 
there is a toric morphism $\p_R: X \ra Y = X(\D^*)$ which is an elementary contraction 
in the sense of Mori theory: $\p_{R^*}(\ccO_X)$ and
$\p_R(C)=point \Leftrightarrow C \in R$. 
Furthermore, let

$$
\begin{array}{cccc}
&A & \ra & B\\
&\cap&&\cap\\
\p_R:&X &\ra&Y\\
\end{array}
$$
be the loci on which $\p_R$ is not an isomorphism;
$A$ and $B$ are the irreducible closed toric strata 
(\ldots), and $\p_R:A\ra B$ 
is a flat morphism, all of whose fibres are 
\textbf{weighted projective spaces}.

\end{theo}

If we use the standard definition of weighted projective spaces then we see that 
the last sentence is false. In the proof Reid describes the fan 
of possible fibres of such contractions and call them weighted projective space. 
We will show that wps's are not the only toric varieties with fan described by Reid.

In this paper we define fake weighted projective spaces as 
a class of spaces defined in Miles Reid proof. 
Next we will see that wps form a subclass of those - the aim of this article is to understand 
how to distinguish them among all fake weighted projective spaces. 
To give an elegant answer we introduce the notion of fundamental group in codimension 1. 
The statement is that a fwps is a wps if and only if it has 
a trivial fundamental group in codimension 1.

We also we introduce the notion of covering in codimension 1 and 
we show that every fwps has a unique universal covering in codimension 1, 
which is a wps.

Using those facts we describe all possible 
algebraic actions of a cyclic group on a $\P^2$, 
which are free away from finite number of points. 
We prove that any such quotient is a fwps.

\section{Notation and definitions}

In this paper by a variety we understand an algebraic variety over the complex numbers. 
In this section we briefly recall basic notions and notation of the toric geometry, 
most of them taken from \cite{fulton}, in particular:

\medskip

\noindent $N$ is the \textbf{lattice of one parameter subgroups of the torus} 
i.e. $\Z^n$ with the standard scalar product.

\noindent A lattice vector $v$ is a \textbf{primitive vector}, if
any integral sub-multiple of $v$ is not in the lattice i.e.: 
$\forall_{ n \in \N} \ \ \ {v \over n} \notin N$.\\

\noindent $M$ is \textbf{the lattice of characters} which is 
a dual lattice to $N$ i.e.~$M = \Hom(N,\Z)$. \\

\noindent $\s$ is a convex, rational, polyhedral cone i.e.~subset 
of $N_{\R} \simeq \R^n$ such that: 
\begin{itemize}
\item there exist lattice vectors $v_1,\ldots,v_k \in N$ such that 
any point $x \in \s$ can be written as $\sum a_i v_i$, 
where $a_i \in \R_{\ge 0}$
\item does not contain any linear subspace, i.e. $\s \cap (-\s)=0$.
\end{itemize}

Such $\s$ we will also denote by $\cone (v_1, \ldots, v_k)$.\\

\noindent $\s^{\lor}$ is the \textbf{dual cone} of the cone $\s$ i.e.  
$\s^{\lor}=\{v \in M_{\R}=N_{\R}^* : 
\forall x \in \s  \hspace{5pt} v(x) \geq 0 \}$.  \\

\noindent $S_{\s}=\s^{\lor} \cap M$ is the \textbf{monoid} 
associated to the cone $\s$ where the addition agrees with addition in $M$.\\

\noindent $A_{\s}$ is the \textbf{algebra} of the cone $\s$ that is 
the polynomial algebra generated by the monomials from semigroup $S_{\s}$. 
In short:  $A_{\s}=\C[S_{\s}]$.\\

\noindent $U_{\s}$ is the \textbf{affine variety} associated with the cone $\s$ that is 
$U_{\s}=\Spec A_{\s}$.\\

\noindent $\D$ denotes a \textbf{fan} in the lattice $N$ which is 
a set of cones in $N_{\R}$ closed under taking and intersecting faces. 
More precisely, $\D$ satisfies the following conditions:
\begin{itemize}
\item if $\s \in \D$, then all its faces are also in $\D$;
\item if cones $\s, \t \in \D$, then $\s \cap \t$ is a face of both $\s$ and $\t$.
\end{itemize}
\textbf{FACT \cite[thm 1.3]{oda}.} If $\t$ is a face of $\s$, then the inclusion 
$\t \hookrightarrow \s$ induces an inclusion 
$U_{\t} \hookrightarrow U_{\s}$ onto an open subset.\\

\noindent A \textbf{toric variety} is a variety defined by a fan $\D$ in the following way: 
$X(\D)=\coprod_{\s \in \D} U_{\s} \slash \sim$, 
where $\sim$ denotes the identification of the open subsets corresponding to common faces. \\

\noindent \textbf{FACT \cite[\S 1.4 p 21 and \S 2.1, p 29]{fulton}.} 
Toric variety defined as above is separable and normal.\\

\noindent Here we introduce some additional, handy notation:\\
For a cone $\s$ by $\D(\s)$ we mean the fan consisting of a cone $\s$ and all its faces.\\
For a fan $\D$ by $\D(1)$ we mean one-dimensional skeleton of the fan $\D$, 
i.e. $\D(1) = \{ \s \in \D : \dim \s = 1\}$.\\
For a collection of vectors $v_1,\ldots,v_n$ by $\conv  \Li v $ we mean
the convex hull of the vectors $v_1,\ldots,v_n$ i.e.~
$\conv  \Li v = \bigl\{ \sum _{i=1}^n t_i \cdot v_i : \sum_{i=1}^n t_i = 1\bigr\}$.\\

\begin{rem}
We will not distinguish between a one dimensional cone $\s$ in $N_{\R}$ 
and the primitive vector generating $\s \cap N$.
\end{rem}

\hfill\\
\noindent Let $\D$ and $\D'$ be fans respectively in $N_{\R}$ and in $N_{\R}'$. \\
\noindent A \textbf{map of fans} is a lattice homomorphism $\p : N' \ra N$, such that
$$
\forall_{\s \in \D}\ \  \exists_{\s' \in \D'}\ \ \  \p (\s) \subset \s'.
$$
We will denote it also by $\p : \D' \ra \D$.
\noindent Recall that any map of fans induces maps of the: 
lattices of characters $M \ra M'$, monoids $S_{\s} \ra S_{\s'}$, 
algebras $A_{\s} \ra A_{\s'}$ and of the affine varieties $U_{\s'} \ra U_{\s}$ 
gluing together to a map  $X(\D') \ra X(\D)$.\\

\noindent A \textbf{toric map} is a map coming from a map of fans.\\

\noindent \textbf{Notation.} If $\p : N' \ra N$, $\p : \D' \ra \D$ is a map of fans, 
then the corresponding toric morphism will also be denoted by $\p$.\\

\noindent \textbf{FACT.}  On a toric variety $X(\D)$ there is an action 
of the big torus $T_{N} \izo N \otimes_{\Z} \C^*$. 
Moreover to any $k$-dimensional cone $\s$ we can associate $n-k$ dimensional orbit 
$\ccO_{\s}\izo ( \C^*)^{n-k}$ (\cite[\S 3.1]{fulton}), 
where the torus $(\s + (-\s)) \cap N \otimes \C^*$ acts trivially.\\

\noindent A \textbf{1-parameter subgroup} 
$\lambda_v : \C^* \ra N \otimes_{\Z} \C^*$ associated with a vector $v \in N$ 
is embedded in the big torus in the following way: 
\mbox{$\lambda_v(t)=(t^{v_1},\ldots,t^{v_n})$}\\

\noindent For any cone $\s \in \D$ a \textbf{toric stratum} $V(\s)$ 
is the closure of the orbit $\ccO_{\s}$.

\section{$\pi_1^1$ - the fundamental group in codim 1}

In this section we introduce the notion of the fundamental group in codimension 1 
and will prove some of its properties.

\begin{defin}
For an irreducible complex algebraic variety $X$ we define 
\textbf{fundamental group in codimension 1} as the inverse limit:
 $$\pi_1^1(X)=\lim_{\longleftarrow \atop  {U \subset X} } \pi_1(U)$$
where $U$ goes over all open, non-empty, algebraic subsets of $X$ 
such that $$\codim_{\C}(X\backslash U) \geq 2$$
\end{defin}

\begin{rem}
We will say $X$ is 1-connected in codimension 1, if \mbox{$\pi_1^1(X)=1$}.
\end{rem}

\begin{rem}
The choice of the base point is not important in our setup, 
because all the varieties under consideration are irreducible 
and hence arc-wise-connected.
\end{rem}
\hfill\\

\noindent The following theorem together with corollary 
\ref{p11_for_smooth} and \ref{p11_for_normal} explain 
how to calculate $\pi_1^1 $ of a given normal variety.

\begin{theo}\label{pi11_ma_sens}
If $V$ is a closed subvariety of a smooth variety $X$ of the complex codimension at least $2$ 
($\codim_{\C} V \geq 2$), then $\pi_1(X)=\pi_1(X\backslash V)$.
\end{theo}

\noindent For simplicity of the proof we assume that $V$ is smooth. 
We will use the following three facts taken from  the book \cite{DNF}. 
These are respectively {\cite[III.2.\S 10 thm 1]{DNF}}, 
{\cite[III.3.\S 12 thm 1]{DNF}}, {\cite[III.2.\S 10 thm 3]{DNF}}:

\begin{theo}[Approximation by a smooth map] \label{aprox}
Let $X$ be a Riemannian manifold with metric $\rho$.
If $Y$ is a compact variety, $f : Y \ra X$ is a continuous map, 
then it can be approximated by a smooth map with given precision i.e. 
$\forall \e>0$ there exists smooth map $g : Y \ra X$ such that $\r(f,g)<\e$.
\end{theo}

\begin{theo}[Maps close to each other are homotopy equivalent]\label{e-homo}
Let $X$ and $Y$ be smooth, compact Riemannian manifolds and $f,g : Y \ra X$ be continuous maps. 
If there exists $\e >0$ such that if $\rho(f,g)< \e$, then $f$ and $g$ are homotopy equivalent.
\end{theo}

\begin{theo}[Transversality]
Let $X$, $Y$ be ${\mathcal\C}^\infty$-smooth manifolds, $V \subset X$ smooth closed submanifold. 
If $f: Y \ra X$ is smooth, then in any neighbourhood of $f$
we can find $g : Y \ra X$  transverse to $V$.
\end{theo}

\begin{prf}{ (of Theorem \ref{pi11_ma_sens}) }
We will show that the inclusion $X \backslash V \hookrightarrow X$ induces isomorphism 
of the fundamental groups. Let us look at a part of the long homotopy exact sequence 
of the pair $(X, X \backslash V)$:
$$
 \pi_2(X,X\backslash V) \ra \pi_1(X \backslash V) \ra \pi_1(X) \ra \pi_1(X,X \backslash V)
$$
We will show that $\pi_2(X,X\backslash V)=\pi_1(X,X\backslash V) = 1$.\\
Let us take a representation of any element from $\pi_2(X,X\backslash V)$ that is 
a continuous map of pairs $\omega : (D^2,S^1) \ra (X, X \backslash V)$. \\
\noindent Note that as $S^1$ is compact and $X \backslash V$ is open, 
the maps close enough to $\omega$ are maps of pairs $(D^2, S^1) \ra (X,X \backslash V)$ 
and are homotopy equivalent with $\omega$ (theorem \ref{e-homo}). 
Using approximation theorem \ref{aprox} we can assume, 
that $\omega$ is smooth. 
Moreover, using the transversality theorem, 
we can assume that $\omega$ is transverse to $V$. 
As the sum of the dimensions over $\R$ of $D^2$ and $V$ is strictly smaller 
than the dimension of $X$, then by  transversality we know that 
the image of $\omega$ does not intersect $V$. 
Thus $\omega : (D^2,S^1) \ra (X \backslash V, X\backslash V)$. 
We have shown that $\omega$ represents the neutral element in $\pi_2(X,X\backslash V)$.\\
With the same technique we can show that $\pi_1(X,X\backslash V)=1$.
\end{prf}

\begin{rem}
In the appendix we have written the definitions and theorems which we need
for the proof of theorem \ref{pi11_ma_sens} without the smoothness assumption for $V$. 
The proof in that more general situation stays the same, 
only in the last step we use more general transversality theorem.
\end{rem}

\noindent We will state some useful corollaries from theorem \ref{pi11_ma_sens}.

\begin{cor}\label{p11_for_smooth}
If $X$ is smooth, then $\pi_1^1(X)=\pi_1(X)$.
\end{cor}
\begin{prf}{}
In \ref{pi11_ma_sens} we have shown that $\pi_1(X)=\pi_1(X\backslash V)$ 
for $V$ of codimension at least $2$. 
From the definition
$$\pi_1^1(X)=\lim_{\longleftarrow \atop  {V \subset X} } 
\pi_1(X \backslash V) = \lim_{\longleftarrow \atop  {V \subset X} } \pi_1(X) = \pi_1(X).$$
\end{prf}

\hfill\\
\noindent \textbf{Fact \cite[II.\S 5 thm 3]{szafarewicz}.} 
Any normal variety $X$ is smooth in codimension~1, 
i.e. $\codim (\Sing(X)) \ge 2$.

\begin{cor}\label{p11_for_normal}
If $X$ is a normal variety, $\Sing(X)$ the set of its singular points, 
$X_0:=X\backslash \Sing(X)$ set of smooth points, 
then $\pi_1^1(X)=\pi_1(X_0)$.
\end{cor}
\begin{prf}{}
From the proof of \ref{pi11_ma_sens} we conclude that 
$\pi_1(X_0)=\pi_1(X_0 \backslash V)$ 
for $V \subset X$ of codimension at least $2$, 
so the inverse limit in the definition of $\pi_1^1$ has a realisation on the set $X_0$.
\end{prf}

\noindent In order to examine the nature of singular points on 
complex normal surfaces Mumford in (\cite{mumford_normal_sings_surf}) 
introduced the local fundamental group. 
He defined it as direct limit over pointed neighbourhoods (in the analytic topology) 
$$
\pi_x= \lim_{\longleftarrow \atop  {U \ni x} } \pi_1(U \backslash \{x\})
$$
In the case of surfaces our fundamental group in codimension 1 
is defined in the same way in the algebraic category i.e.~
all the sets $U$ are open in the Zariski topology.

Mumford proved that rational surface singularities are characterised 
by their local fundamental groups. More precisely:

\begin{defin}
We call $(X,x)$ a \textbf{rational singularity} if for some resolution of singularities 
$\pi : X' \ra X$ higher derivatives of the direct image functor in $\pi$ vanish: 
$R^i\pi_*(\ccO_{X'})=0$ for $i>0$.
\end{defin}

\noindent Rationality of $(X,x)$ can be characterised as follows: 
in some resolution of singularities $(\hat X,\bigcup E_i) \ra (X,x)$ 
exceptional curves $E_i$ are smooth, rational, intersect transversely 
and their intersection graph is acyclic.

\begin{theo}[Mumford]
If (X,x) is a normal rational complex singularity and $\pi_x(X)=0$, then $x$ is a smooth point.
\end{theo}

\noindent Another similarity with Mumford's approach
is that in this paper we deal with toric varieties, 
which have only rational singularities.
\noindent We will see that any wps is 1-connected in codimension 1. 
A wps surface is usually singular in 3 points (\cite{dimca}). 
It means that on a wps any loop around one singular point 
can be contracted without touching singularities.

\begin{defin}
A finite, surjective morphism $\p : Y \ra X$ is called 
\textbf{covering in codimension 1} 
if it is unramified in codimension 1. 
More precisely, there exists a subvariety $V \subset X$ such, that:
\begin{itemize}
\item $\codim(V) \ge 2$
\item$\p\mid_{Y_1} : Y_1 \ra X\backslash V$ is a topological cover, 
where $Y_1=\p^{-1}(X\backslash V)$.
\end{itemize}
\end{defin}

\begin{rem}
Let $Y \ra X$ be a covering in codimension 1. 
It is universal (in the usual, categorical sense) if and only if 
$Y$ is one-connected in codimension 1, i.e. $\pi_1^1(Y)=0$.
\end{rem}

\begin{lem}
Let $p : \tilde X_0 \ra X_0$ be the topological universal covering 
of the set of smooth points $X_0$ of a normal variety $X$ such, that 
$p$ and $\tilde X_0$ are algebraic.
Then $p$ extends to a universal covering in codimension 1.
\end{lem}

\begin{prf}{}
Let us first observe that it is enough to solve the problem for an affine $X$ 
(compare \cite[\S 2.14]{iitaka}). If $X$ satisfies the above conditions, 
we have following diagrams of varieties and of algebras:
$$
\begin{array}{ccccccc}
&&X&\hfill&         &&A     \\
&&\uparrow&\hspace{40 pt}&      &&\downarrow    \\
\tilde X_0 & \stackrel p \ra & X_0&\hfill&  \tilde A_0&\stackrel p \la &A_0\\
\end{array}
$$
We want to find $\tilde X$ and $\tilde A$, extending these diagrams. 
Denote by $\tilde K$ the quotient field of $\tilde A_0$, 
by $\tilde A$ integral closure of the ring $p(A)$ in $\tilde K$ 
and by $\tilde X = \Spec \tilde A$. 
Then inclusion $A \hookrightarrow \tilde A$ gives the required extension.
\end{prf}

\section{The group $\pi_1^1$ for toric varieties}
\noindent In this section we calculate the fundamental group in codimension 1 
of a toric variety in terms of its combinatorial description. 
See \cite{fulton} for standard toric geometry used in this paragraph.

\begin{lem}\label{pi1U_na_pi1X}
Let $X$ be a normal, irreducible variety, $U$ its open, dense subset in the Zariski topology. 
Then the inclusion $U \hookrightarrow X$ induces surjection of the fundamental groups: 
$\pi_1(U) \twoheadrightarrow \pi_1(X)$,
\end{lem}

\begin{rem}
This lemma generalises the part of theorem \ref{pi11_ma_sens} 
where we have shown that $\pi_1(X,U)$ 
appearing in the following long exact homotopy sequence of the pair 
$(X,U)$ is trivial:
$$\pi_1(U) \twoheadrightarrow \pi_1(X) \ra \pi_1(X,U) \ra \pi_0(U) = 1$$
\end{rem}

\begin{prf}{(of lemma \ref{pi1U_na_pi1X})}
 Let $p : Y \ra X$ be the universal covering. 
Let us restrict it to a Zariski open set $U$, 
i.e.~let us consider $p_1 : Y_1 \ra U$, where $Y_1 =p^{-1}(U)$. 
Our variety $X$ is normal and irreducible, so it is locally irreducible as an analytic space, 
i.e.~cutting out subset of a strictly smaller dimension (over $\C$) 
does not disconnect the connected open subsets. 
Thus, as $p$ is an open map, $Y_1$ is connected. 
Let us consider long exact homotopy sequence for fibrations \cite[\S 7, \S 2, Thm 10]{spanier}
- write it for $p$ and $p_1$:
$$
\begin{array}{ccc ccc ccc c}
 & \pi_1(Y_1) & \stackrel{p_1^*}\ra & \pi_1(U)   & \ra & \pi_0(F)   & \ra &  \pi_0(Y_1)  & =  & 1\\
   & \downarrow && \downarrow   &     & \Arrowvert &     &  \Arrowvert && \\
 1 =  & \pi_1(Y)   & \stackrel{p^*}\ra & \pi_1(X) & \ra & \pi_0(F)   & \ra &  \pi_0(Y)   & = & 1 \\
\end{array}
$$
By a usual diagram chase we conclude that the homomorphism in question is surjective.
\end{prf}

\begin{defin}
Let $v$ be a vector in the lattice $N$, let $\lambda_v : \C^* \ra T_N$ be the corresponding 
one-parameter subgroup. With the vector $v$ there is an 
\textbf{associated loop} in the torus $T_N$, which we will also denote by $\lambda_v$. 
It is given by the inclusion \mbox{$S^1  \izo U(1) \hookrightarrow \C^*$}, 
associated with one-parameter subgroup $\lambda_v$ via 
\mbox{$U(1)=\{t : |t| =1 \}$}.
\end{defin}

\begin{cor}\label{petle_alg}
Let $\D$ be a fan in $N_{\R}$. Any loop in $\pi_1(X(\D))$ is homotopy equivalent 
with an algebraic 1-cycle, i.e.~associated to a vector $v \in N$.
\end{cor}

\begin{prf}{}
Any toric variety is normal, irreducible and 
contains as an open dense subset the big torus. 
Thus by theorem \ref{pi1U_na_pi1X} we know that 
$\pi_1(T_N) \twoheadrightarrow  \pi_1(X(\D))$, 
and $\pi_1(T_N) \izo N$ is of course generated by 
$\lambda_{e_1}, \ldots, \lambda_{e_n}$.
\end{prf}

\begin{lem}\label{pi11_tor}
For a fan $\D$ and toric variety $X(\D)$ we have
$$\pi_1^1(X(\D))=\pi_1(X(\D(1)))$$.
\end{lem}
\begin{prf}{}
$X(\D(1))$ is an open subset of $X(\D)$ whose complement is of codimension at least 2: 
$\codim_{\C}(X(\D) \backslash X(\D(1)))\ge 2$.
So $$\pi_1^1(X(\D))=\pi_1(X(\D(1)))$$ by Theorem \ref{pi11_ma_sens}.
Moreover $X(\D(1))$ is smooth so by Corollary \ref{p11_for_smooth} the claim follows.
\end{prf}

\begin{lem}\label{petle_sciagalne}
If a vector $v$ lies in the relative interior of a cone  $\s \in \D$, 
then the loop $\lambda_v$ is contractible in $U_\s$.
\end{lem}

\begin{prf}
If $v$ lies in the relative interior of a cone  $\s \in \D$, 
then there exists the limit  
$\lim_{t \ra 0} \lambda_v(t) = x_{\s}$ (\cite[\S 2.3]{fulton}). 
Thus the map $\lambda_v : \C^* \ra T_N \mono U_\s$ can be extended 
to a map  $\widetilde {\lambda_v} : \C \ra U_\s$, 
which implies that the loop $\lambda_v$ is contractible.
\end{prf}

\begin{lem}\label{pi1_x1}
Let  $N$ be a lattice of rank $n$, let $\D$ be any fan in  $N_{\R}$, 
with rays generated by primitive vectors $v_1, \ldots, v_k$. Then:
$$\pi_1 \Bigl( X \bigl( \D(1) \bigr) \Bigr)=N \slash N'$$
where $\D(1)$ is the one-dimensional skeleton of $\D$, 
and $N'\subset N$ sublattice generated by  $v_1, \ldots, v_k$.
\end{lem}

\begin{prf}{}
From the corollary \ref{petle_alg} we know that $\pi_1(X(\D(1))) $ is generated by 
$\{ \lambda_v : v \in N \}$. Also, from lemma \ref{petle_sciagalne}, 
we know that the loops $\lambda_{v_1}, \ldots , \lambda_{v_k} $ are contractible. 
Moreover, for any vector $v_i$ we have an isomorphism 
$$\C^* \times \ldots \times \C^* \times \C \izo X(\cone(v_i))$$
thus $\pi_1(X(\cone(v_i))) = N / \Z \cdot v_i$.
Now $X(\D(1))$ is covered by open sets $X(\cone(v_i))$, 
of which fundamental groups we already know. 
Intersection of any two of those sets is the big torus. 
We use the van Kampen theorem:
$$
\pi_1 \Bigl( X \bigl( \D(1) \bigr) \Bigr)
=\Biggl({\prod  \pi_1\Bigl(X\bigl(\cone(v_i)\bigr)\Bigr)}\Biggr) \big\slash_\sim  
= \Bigl(\prod \bigl(N \slash \Z \cdot v_i\bigr)\Bigr) \big\slash_\sim 
= N \slash N'$$
where $\sim$ denotes identifications along the map coming from 
$T_N \hookrightarrow X(\cone(v_i))$, that is coming from the surjection: 
$N \twoheadrightarrow N \slash {\Z\cdot v_i}$.
\end{prf}

\noindent We summarise the results of this section 
in the following theorem characterising 
fundamental group in codimension 1 for toric varieties:

\begin{theo}[$\pi_1^1$ for toric varieties]\label{pi11tor}
Let $\D$ be a fan in a lattice $N$. Then:
$$\pi_1^1\bigl(X(\D)\bigr) = N \slash N' $$
where $N'$ is the sublattice of $N$ generated by vectors from 1-dimensional cones in~$\D$.
\end{theo}

\section{Weighted projective spaces}

In this section we will recall the definition of a weighted projective space. 
Later we will calculate its fundamental group in codimension~1.
\begin{notation}\label{not}
Let $v_0, \ldots, v_n$ be primitive vectors in a lattice $N$ of  rank $n$ such, that $0 \in \operatorname{int} (\conv \Lo v)$. By $\D \Lo v$ we denote the fan consisting of the cones $\s_i=\cone (v_0,\stackrel {\stackrel i \lor } \ldots, v_n)$ and all their faces.
\end{notation}

\begin{defin}
Let $a_0,\ldots,a_n$ be positive integral numbers. 
A \textbf{weighted projective space} with weights $(a_0,\ldots,a_n)$ 
is the orbit space of the following action of $\C^*$ on $\C^{n+1} \backslash \0$:
$$t \cdot (x_0, \ldots ,x_n) = (t^ {a_0} \cdot x_0,\ldots,t^{a_n} \cdot x_n)$$ 
and it is denoted by $\P(a_0, \ldots, a_n)$. 
We will denote the class of weighted projective spaces by $w\P$.
\end{defin}

\begin{rem}
For any natural number $k$ there is an isomorphism:
$$\P \Lo {ka} \izo \P \Lo a.$$
therefore we can assume that $GCD \Lo a = 1$, 
but even with this assumption the numbers $a_0,\ldots,a_n$ 
are not uniquely defined (\cite{fletcher}, \cite{dimca}, \cite{dolgachev}).
\end{rem}

\noindent The following proposition describes the structure
of a toric variety on a w.p.s.

\begin{prop}\label{wP_alg}
Let $N$ be a lattice in $\R^n$ generated by primitive vectors $v_0,\ldots,v_n$ 
such, that  $a_0v_0+\ldots+a_nv_n=0$. Then 
$$X \bigl(\D \Lo v \bigr) \izo \P \Lo a$$
\end{prop}
\begin{prf}{}
Denote by $\D$ the fan $\D \Lo v$, 
by $p : \C^{n+1} \backslash \0 \ra \P \Lo a$ 
the quotient map from the above definition. 
The aim is to prove that, $p$ is toric. 
First we explain what the fan of
$\C^{n+1} \backslash \0$ looks like.
Let $N'$ be a lattice of rank $n+1$ in $\R^{n+1}$, 
with the standard basis $e_0,\ldots,e_n$. 
Let $\D'$ be the fan constructed in the same way as $\D$,
i.e. $\D'$ consists of the cones
$\{\cone (e_0,\stackrel {\stackrel i \lor} \ldots, e_n)\}$. 
Then $\D'$ is the fan of $\C^{n+1} \backslash \0 $.
We define a map
$\p : \C^{n+1} \backslash \0 \ra X(\D)$ on the level of lattices by the formula
$\p(e_i)=v_i$. We will show that this coincides with $p$. 
It is straightforwards that $\p$ is the projection along the vector $\Lo a$, 
so its kernel is generated by this vector and 
we have the following short exact sequence:
$$
0 \ra \Z \cdot \Lo a  \ra N' \ra N \ra 0.
$$
Now we apply the functor $\otimes_{\Z} \C^*$ to this sequence in order to obtain
$\p : T_N \ra T_{N'}$ on the level of tori.
Because  $\p : N' \ra N$ is surjective we get the exact sequence:
$$
0 \ra \bigl(\Z \cdot \Lo a  \otimes_{\Z} \C^*  \izo \C^*\bigr) \ra 
N' \otimes_{\Z} \C^* \ra N \otimes_{\Z} \C^* \ra 0
$$
where $\C^* \hookrightarrow N'\otimes_{\Z} \C^*$ 
is given by: $t \mapsto  \Lo a \otimes t$. 
Thus the kernel of $\p$ as a map of tori is the one-parameter subgroup
corresponding to the vector $\Lo a$.
We have shown that for the big torus we are dividing out by the $\C^*$-action, 
and hence $p \mid_{T_{N'}} = \p \mid_{T_{N'}}$, which means $\p=p$. 
\end{prf}

\noindent We note, that any weighted projective space
$\P \Lo a$ is obtained as a quotient of $\P^n$ by an action of some finite group.
This map is usually branched in codimension 1.
\begin{prop}\label{wp}
Any weighted projective space is the orbit space of an action of a finite group 
on a projective space of the same dimension:
if \mbox{$GCD(a_0,\ldots,a_n)=1$}, then $\P(a_0,\ldots,a_n)=\P^n \slash G$, 
where $G=\Z_{a_0}\oplus \ldots \oplus \Z_{a_n}$ acts coordinate-wise (see equation below).
\end{prop}
\begin{prf}{}
We will construct toric maps $q$ and $q'$ in the following diagram:

$$
\begin{array}{ccc}
\C^{n+1} \backslash \0 & \stackrel {q'} \mr & \C^{n+1} \backslash \0 \\
\mapdown {\p_1} && \mapdown \p\\
\P^n & \stackrel q \mr & \P \Lo a \\
\end{array}
$$
where $\p$ is a projection defined as in the proof of proposition \ref{wP_alg} and 
$\p_1$ is the usual dividing out by scalar multiplication.\\
On lattices $\p : N' \ra N$ (respectively $\p_1 : N_1' \ra N_1$)
is a projection  along the vector $\Lo a$ 
(respectively $(1,\ldots,1)$). 
We will choose a basis for the sublattice $N_1'$ 
of the lattice $N'$ in such a way that, the vector $(1,\ldots,1)$ 
from the lattice $N_1'$ will become the vector $\Lo a$ 
from the lattice $N'$. \\
Let $N_1'=\Z a_0  e_0 \oplus \ldots \oplus \Z \ a_n  e_n$ 
and let $q' : N_1' \hookrightarrow N'$ be the inclusion of lattices. 
On the level of varieties $q'$ corresponds to a quotient by the action of 
$G=\Z_{a_0}\oplus \ldots \oplus \Z_{a_n}$: 
let $\e_i$ be the primitive root of $1$ of degree $a_i$ 
and at the same time generator of $\Z_{a_i} \subset G$.
 Then the action is as follows:
$$
(1,\ldots,\e_i,\ldots,1) \Lo x = (x_0,\ldots,\e_i \cdot x_i,\ldots,x_n)
$$
and it descends to the required action on $\P^n$.
Thus the earlier described diagram of lattices can be completed to a commutative one:
$$
\begin{array}{ccc}
 N_1'  & \stackrel {q'}\mr & N' \\
 \mapdown {\p_1}  && \mapdown \p\\
N_1 & \stackrel q \mr & N \\
\end{array}
$$
In that way we have constructed the map $q : \P^n \ra \P \Lo a$ completing the initial diagram.
\end{prf}

\begin{rem}
The numbers $a_0,\ldots,a_n$ are not uniquely defined, 
nor are the maps in question.
\end{rem}
\hfill\\
\noindent The last theorem of this section will summarise results 
of the previous two sections in the useful form:
\begin{theo}\label{pi11wp}
Weighted projective space $\P(a_0,\ldots,a_n)$ is one-connected in codimension 1
i.e. $\pi_1^1(\P \Lo a)=0$.
\end{theo}

\begin{prf}{} 
This follows from the theorem \ref{pi11tor}, 
which says that $\pi_1^1(X(\D))=N \slash N'$. 
For  $\P \Lo a$ lattices $N$ and $N'$ are the same ---
 $N'$ is defined as sublattice $N$ generated by the vectors 
$v_0,\ldots,v_n \in \D(1)$.
Proposition \ref{wP_alg} says that the vectors $v_0,\ldots,v_n$ 
generate the lattice $N$.
\end{prf}

\section{Fake weighted projective spaces}

In this section we define fake weighted projective spaces (fwps).
The fan of a fake weighted projective space is described 
in the Reid's proof of theorem \ref{Reid},
and it is almost the same as a fan  of a weighted projective space.

We already have introduced notation for those (\ref{not}),
when we have been defining weighted projective spaces.
The only difference between these two classes of varieties, 
when we look at fans and lattices is
that in fwps case the vectors defining the fan need not to generate the lattice $N$.

\begin{defin}
We call $X(\D \Lo v)$ a \textbf{fake weighted projective space} for any vectors
satisfying assumptions of the notation \ref{not}. 
We denote the class of fake weighted projective spaces by $t\P$, 
or by $t\P^n$ if we mean the subclass of varieties of dimension n in $t\P$.
\end{defin}

\begin{ex}\label{ex}
The class $t\P$ is significantly bigger than the class $w\P$:
let us take the fan generated by the following vectors in the 
standard $N=\Z^2$ lattice:
$$v_0=(1,-1), \ \ v_1=(1,2), \ \ v_2=(-2,1)$$
The generated lattice has index 3, thus $X(\D(v_0,v_1,v_2))$ 
is not a wps, but of course is in $t\P$. 
As we will see later it is a $\P^2$ divided by the following action of $\Z_3$ 
$$\e \cdot (x_0,\ x_1,\ x_2) = (x_0,\  \e \cdot x_1,\ \e^2 \cdot x_1)$$
where $\e$ is primitive root of unity of degree 3.
\end{ex}

\noindent Before we state the main result of this section 
we need to prove the following lemma:

\begin{lem}\label{inclusion}
Let $\hat N \subset N$ be both of rank $n$.
Let $\hat \D$ and $\D$ the same fan embedded respectively in
the lattice $\hat N$ and $N$. 
Then the finite group $G=N \slash \hat N$ acts on $X(\hat \D)$ 
and the quotient map $\xi : X(\hat \D) \ra X(\D)$ 
comes from the inclusion of lattices $\xi : \hat N \ra N$.
\end{lem}
\begin{prf}{}
The following sequence is exact:
$$0 \ra \hat N \ra N \ra G \ra 0.$$
We apply the functor $\otimes_{\Z}\C^*$ to it and we get:
$$0 \ra Tor^1_{\Z}(G ,\C^*) \ra T_{\hat N} \ra T_N \ra G \otimes_\Z \C^* \ra 0.$$
We can rewrite the above sequence as:
$$0 \ra G  \ra T_{\hat N} \ra T_N \ra 0.$$
Thus on the tori we have a topological covering $\xi : T_{\hat N} \ra T_N$ 
with Galois group $G$. The group $G$ acts on $T_{\hat N}$  by left multiplication.
This action can be extended to $X(\hat \D)$. 
We will show, that $\xi$ is the quotient map.
Note that it is enough to find a covering by affine subsets with this property.\\
\indent Fix $\s \in \D$. 
We will show, that the restriction $\xi \mid_{U_{\hat \s}} : U_{\hat \s} \ra U_{\s}$  
is the quotient map.
The action of $G$ on the torus $T_{\hat N}$ induces the action of $G$ 
on the function algebra of the torus $A_{\hat 0}$ 
and the map $\xi : A_0 \ra A_{\hat 0}$ 
is an inclusion on the $G$ - invariant subalgebra. 
As $G$ is a subgroup of the torus, it acts also on the algebra $A_{\hat \s}$, 
thus the restriction $\xi \mid_{A_{\s}} : A_{\s} \ra A_{\hat \s}$ 
is an inclusion on the $G$-invariant subalgebra as well.
We conclude that $\xi : X(\hat \D) \ra X(\D)$ is the quotient map.
\end{prf}

\begin{theo}\label{uni}
For any fake weighted projective space $X(\D \Lo v)$ 
there exist a unique universal covering in codimension 1,
which is a weighted projective space $\xi : \P \Lo a \ra X(\D \Lo v)$. 
Moreover the numbers $a_0,\ldots,a_n$ are positive, integral
and they satisfy the condition: $\sum_{i=0}^n a_iv_i = 0$. 
\end{theo}

\begin{prf}{}
The covering is induced by lattice inclusion.
Let $\hat N$ be the sublattice of $N$ generated by vectors $v_0,\ldots,v_n$. 
These vectors are dependent in $N$, moreover $0 \in \conv \Lo v$.
Let $a_0,\ldots,a_n$ be positive integers such, that  $\sum_{i=0}^n a_iv_i = 0$. 
Let us consider the fan $\D \Lo v$ in the lattices $\hat N$ and $N$. 
Denote it respectively by $\hat \D$ and $\D$. 
Of course $X(\hat \D)$ is the wps $\P \Lo a$.  
The lemma \ref{inclusion} says that, the map $\xi : X(\hat \D) \ra X(\D)$ 
induced by the inclusion of lattices is a quotient map 
and that there is an induced group action of $N \slash \hat N$ on $\P \Lo a$.
We need to show that $\xi$ is unbranched in codimension 1.
We know that:

\begin{itemize}
\item branching points of $\xi$ are in 1-1 correspondence with
the  points with non-trivial isotropy group of the action of $N \slash \hat N$
\item all points in an orbit of the action of big torus have the same isotropy group
\item any orbit $\ccO_{\t}$ is contained in the affine set $U_{\t}$.
\end{itemize}

Let us fix a $1$-dimensional cone $\t$ (in other words a vector $v \in \{v_0,\ldots,v_n\}$) 
and let us check that $\xi \mid_{U_{\hat \t}} : U_{\hat \t}\ra U_{\t}$ 
is unbranched on $\ccO_{\t}$ 
(where $\hat \t$ denotes the cone $\t$ as a subset of $\hat N$). 
The vector $v$ is primitive both in the lattice $N$ and in the lattice $\hat N$. 
Thus $U_{\hat \t}\izo U_{\t}\izo \C \times (\C^*)^{n-1}$ 
and also $\ccO_{\hat \t}\izo\ccO_{\t}\izo (\C^*)^{n-1}$. 
Moreover the coordinates can be chosen in such a way that
the map $\xi : U_{\hat \t} \ra U_{\t}$ becomes the product
of identity and the already described map on tori.
Consider the following commutative diagram of short exact sequences:
$$
\begin{array}{ccccc ccccc}
0 & \ra & \Z \cdot v & \ra & \hat N & \ra & \hat N \slash \Z \cdot v & \ra & 0\\
  &      & \Arrowvert &    & \downarrow & &\Arrowvert & & \\
0 & \ra & \Z \cdot v & \ra &  N & \ra & N \slash \Z \cdot v & \ra & 0\\
\end{array}
$$
If we split the top sequence and the lower one accordingly,
then the splittings will match. 
Therefore we get the following commutative diagram:
$$
\begin{array}{ccccc}
\C 	   & \ra & U_{\hat \t} & \ra & \ccO_{\hat \t} \\
\Arrowvert &     & \downarrow  &     & \downarrow  \\
\C 	   & \ra & U_{\t}      & \ra & \ccO_{\t}\\
\end{array}
$$
where $\ccO_{\hat \t} \ra \ccO_{\t}$ is a topological cover
with the  Galois group $N \slash \hat N$.
\noindent Finally $\xi$ is a covering in codimension 1 
and by theorem \ref{pi11wp} it is a universal covering in codimension 1.
\end{prf}

\begin{theo}[Characterisasion of the class $w\P$ in the class $t\P$]
If $X$ is a fake weighted projective space then it is a weighted projective space 
if and only if it is 1-connected in codimension 1. In short:
$$X \in t\P \ \ \Longrightarrow \ \ (X \in w\P \Longleftrightarrow \pi_1^1(X)=1)$$.
\end{theo}
\begin{prf}{}
By theorem \ref{uni} we have a uniquely defined universal covering by
a weighted projective space, which is trivial exactly for the class of wps.
\end{prf}

\begin{cor}
We can summarise the results of \ref{uni} and of \ref{wp} by saying that
for any fake weighted projective space $X$ there exists a finite map $\psi : \P^n \ra X$. 
This is not uniquely defined.
\end{cor}

\noindent Before we proceed, let us note that fake weighted projective spaces
can actually appear as fibers of extremal contractions, as described by Reid.
Let us take any 2-dimensional fwps (example \ref{ex}). 
From the cone theorem we know that extremal rays in  $NE(X)$ correspond to
toric strata of dimension 1, that is to cones of dimension 1.
Let us take any cone $\s$ of dimension 1.
The corresponding map contracts the curves numerically equivalent with $V(\s)$. 
Let us recall (\cite[\S3.4, p 63]{fulton}) ) the well-known formula 
for the rank of the Picard group for toric variety with simplicial fan:
$\rk \Pic(X(\D)) = \#\D(1) - \rk N$. In our case:
$X \in t\P^2 \Rightarrow \rk \Pic X = 3 - 2 = 1$, 
so the target of our contraction is a point.

\section{Application - quotients of  $\P^2$ by a cyclic group}

To illustrate  theorem \ref{uni}, we will find explicitly 
all possible actions of a cyclic group on the surface $\P^2$, 
that give the universal covering in codimension~1. 
We will investigate a particular toric map coming from 
the inclusion of lattices as in lemma \ref{inclusion}).

\begin{lem}[Toric quotient of $\C^{n+1}$ by $\Z_r$]\label{gr_sk}
Let $\e$ be a primitive root of unity of degree $r$, 
let $a_0,\ldots,a_n$ be non-negative numbers.
Assume that \mbox{$GCD \Lo a = 1$.}
Fix the action of cyclic group $\Z_r = \langle \e \rangle$ on $\C^{n+1}$ by:
$$\e \cdot (x_0,\ldots, x_n) = (\e^{a_0} \cdot x_0, \ldots, \e^{a_n} \cdot x_n).$$
Then the quotient map $\p : \C^{n+1} \ra  \C^{n+1} \slash \Z_r$ 
is a toric map corresponding to the inclusion of lattices:
$$\xi \ : \ \Z \cdot e_0\oplus \ldots \oplus \Z \cdot e_n 
\ \hookrightarrow \ 
\Z \cdot v + (\Z \cdot e_0 \oplus \ldots \oplus \Z \cdot e_n),$$
where $v={1 \over r}( a_0 \cdot e_0 + \ldots + a_n \cdot e_n)$.
On the level of cones $\xi$ is the identity.
\end{lem}

\begin{prf}{}
Set $\hat N = \Z \cdot e_0 \oplus \ldots \oplus \Z \cdot e_n$ 
and $N=\Z \cdot v +\Z \cdot e_0\oplus \ldots \oplus \Z \cdot e_n$. 
We have the exact sequence: 
$$
0 \ra \hat N \stackrel \xi \ra N \ra \Z_r \ra 0
.$$ 
From lemma \ref{inclusion} we know that $\xi : \C^{n+1} \ra \C^{n+1} \slash \Z_r$ 
is a quotient map and its kernel on the level of tori is $\Z_r$. 
We need to check how this $\Z_r$ is embedded in $T_{\hat N}$.
On tori we have:
$$
0 \ra \Z_r \ra T_{\hat N} \stackrel \xi \ra T_N \ra 0
,$$ 
where $\xi(e_i \otimes t)=e_i  \otimes t$. 
Recall that the isomorphism $T_N \izo N \otimes \C^*$ 
is given by $\Lo v \otimes t = (t^{v_0},\ldots,t^{v_n})$. 
The group generated by 
$(\e^{a_0},\ldots,\e^{a_n})$ is in the kernel:
$$
\xi(\sum_{i=0}^n a_i e_i \otimes \e) = 
\sum_{i=0}^n a_i e_i \otimes \e = r \cdot v \otimes \e = 
v \otimes 1 = 0
$$
(for the  module $N\otimes \C^*$ we use the additive notation). 
We know that $\Z_r$ is the kernel, thus it is generated by
$(\e^{a_0},\ldots,\e^{a_n})$ and the action is diagonal, as required.

\end{prf}

\begin{rem}
In the above proof vector $v$ is not unique:
by taking any $v_1$  determining the same lattice,
we get the same action. In other words if:
$$
v - v_1 \in \hat N
,$$
then the lattices:
$$\Z \cdot v_1 + (\Z \cdot e_0 \oplus \ldots \oplus \Z \cdot e_n) =
\Z \cdot v + (\Z \cdot e_0 \oplus \ldots \oplus \Z \cdot   e_n)$$
are identical, thus the vectors  $v$ and $v_1$ determine the same map.
\end{rem}

\begin{lem}[When is $t\P^2$ covered by $\P^2$ in codimension 1]\label{tp2_p2}
A surface $t\P^2$ is covered by $\P^2$ in codimension 1 if and only if
it is isomorphic to the orbit space of $\Z_r$ acting on $\P^2$ in the following way:
$$
\e \cdot (z_0 : z_1 : z_2) = (z_0 : \e^{a+1} \cdot z_1 : \e^a \cdot z_2)
.$$
In the above we assume $\e$ to be a primitive root of unity of degree r 
and 
$$
GCD(a,r)=GCD(a+1, r) =1
.$$
\end{lem}

\begin{prf}{}
Let $X(\D) \in t\P^2$ and  let $v_0, v_1, v_2$ be the vectors generating $\D(1)$ 
in the lattice $\Z \cdot f_1 \oplus \Z \cdot f_2$. 
We can assume that $v_0=(1,0)$ by applying a lattice automorphism if necessary.
As the covering wps of $X(\D)$ is $\P^2 \izo \P(1,1,1)$, 
the sum of the vectors $v_0,v_1,v_2$ need to be zero: 
$(1,0) + (a,r) + (c,d) = (0,0)$. We get: $1+a+c=0$ and $r+d=0$. 
Thus the fan $\D$ is generated by the vectors $(1,0), (a,r), (-1-a,-r)$. 
These vectors generate lattice of index $r$ in the initial lattice.
Thus on tori we are dividing out by the cyclic group $\Z_r$. 
Now it is enough to see what this action looks like 
in homogeneous coordinates $\P^2$.
We will find an action of $\Z_r$ on $\C^3 \backslash \{0\}$ which descends to $\P^2$. 
Let us write the following  diagram:    
$$
\begin{array}{ccc}
 \C^3 \backslash \{0\} 		& \stackrel {\slash_{\Z_r}} {\mr} 
& \C^3 \backslash \{0\}  \slash \Z_r\\
\mapdown{\slash_{\C^*}} 	& 	& \mapdown{\slash_{\C^*}}\\
\P^2 			& \stackrel {\slash_{\Z_r}} {\mr} & X(\D)  \\
\end{array}
$$
corresponding to the following diagram of lattices:
$$
\begin{array}{cccc}
 \Z e_0 \oplus\Z e_1 \oplus\Z e_2  & \mr  & \Z w + (\Z e_0 \oplus \Z e_1 \oplus\Z e_2)  \\

  \mapdown \slash_{\ (1,1,1)}  & & \mapdown  \slash_{\ (1,1,1)} \\

{(\Z v_0 \oplus \Z v_1 \oplus\Z  v_2) \slash \Z (v_0+v_1+v_2)}  &  \mr  &

{\bigl(\Z  f_2 + (\Z v_0 \oplus \Z  v_1 \oplus \Z  v_2)\bigr) \slash \Z (v_0+v_1+v_2)} \\

\end{array}
,$$
where $w={1 \over r} (a_0,a_1,a_2)$. We already know that:
$$
f_2 = {1 \over r} \bigl((a+1) \cdot v_0 + a \cdot v_2 \bigr)
$$
so we can assume $a_0=0$ and thus $a_1=a+1$, $a_2=a$. 
From the lemma \ref{gr_sk} the action of  $\Z_r$ on $\C^3 \backslash \0$ 
in the homogeneous coordinates is given by the required formula.
\end{prf}

\begin{rem}
Of course the action can be written in various ways - 
we can multiply by a number and choose various primitive roots.
This ambiguity corresponds to that in the proof  \ref{gr_sk}. 
Let illustrate this phenomena on an example. Let $\xi = \sqrt[7]1$. 
Consider the action
$$ 
\xi \cdot (x_0 : x_1 : x_2) = (x_0 : \xi^3 x_1 : \xi^5 x_2)
$$
By substituting $\e = \xi ^2$ we get $\xi^3=\e^5$ and $\xi^5=\e^6$, so the action is:
$$ \e \cdot (x_0 : x_1 : x_2) = (x_0 : \e^5 x_1 : \e^6 x_2)$$
This corresponds to the following ambiguity on the vectors:
$${2 \over 7} (0,5,6) - {1 \over 7} (0,3,5) = (0,1,1)$$
However, if the  $\Z_r$-action is free outside of a finite set, 
it can be written as in lemma \ref{tp2_p2} as follows:
Suppose we have the following action:
$$ \xi \cdot (x_0 : x_1 : x_2) = (\xi^b x_0 : \xi^c x_1 :  x_2) = (x_0 : \xi^{c-b} x_1 : \xi^{-b} x_2)$$
where $GCD(b,r)=GCD(c,r)=GCD(b-c,r)=1$. It is equivalent to the action:
$$ \e \cdot (x_0 : x_1 : x_2) = ( x_0 : \e^{a+1} x_1 :  \e^a x_2)$$
where $\e=\xi^c$, and    $ac = -b  \mod \ r $.
\end{rem}
 
\begin{lem}\textbf{(Lemma \ref{tp2_p2} from different point of view.}\label{lem2})
If an action of $\Z_r$ on $\P^2$ is algebraic and free out of finite number of point, 
then it can be written in some basis in the form:
$$\e \cdot (z_0 : z_1 : z_2) = (z_0 : \e^a \cdot z_1 : \e^{a+1} \cdot z_2)$$
where $\e$ is a root of unity of degree r and $a$ and $a+1$ are co-prime to $r$.
\end{lem}
 
\begin{prf}{}
All the algebraic automorphisms of $\P^n$ form the group $PGL(n+1)$ \cite[Exercise 7.2]{hartshorne}. 
In our case $\Z_r$ is a subgroup of  $PGL(3) = SL(3) \slash \Z_3$. 
Matrix of the generator of $\Z_r$, seen as an element of $SL(3)$,
can be diagonalised because it has finite order.\\
By choosing the coordinates diagonalising the action and using the lemma \ref{tp2_p2} 
we get the required formula.
\end{prf}

\begin{cor}
If a normal surface $S$ fulfils the assumptions of \cite [Thm 2.3]{WW}
and $\pi_1^1(S)$ is a cyclic group, then $S$ is isomorphic 
with a surface of type $t\P^2$ from the lemma \ref{tp2_p2}.
\end{cor}
\newpage
\appendix
\section{Appendix - stratifications}

Here we quote the definitions and theorems needed for 
the proof of theorem \ref{pi11_ma_sens} in its stronger version.
It says that  for smooth  $X$ and not necessarily smooth 
algebraic subvariety $V$ of codimension at least 2, 
one has $\pi_1(X \backslash V) = \pi_1(X)$. 

\begin{defin}
Let $V$ be a closed subvariety of a smooth variety  $X$. 
We say that there exists a \textbf{Whitney stratification of $V$} 
if $V$ can be written as a sum $V=\bigcup_{i \in I} S_i$, 
where $I$ is a partially ordered set and the following conditions hold:

\item{-} $S_i$ is  locally closed, smooth subvariety of $X$,
\item{-} if $\alpha < \beta$ the Whitney conditions A and B
for the pair $S_{\alpha}$, $S_{\beta}$ hold:
for any sequence   $x_i \in S_{\beta}$ 
convergent to the point $y \in S_{\alpha}$ and
for any sequence of points $y_i \in S_{\alpha}$ 
convergent to the same point $y \in S_{\alpha}$
let us denote secants by $l_i=x_i y_i$. 
Suppose that $l_i$ converge to the line $l$ 
and the tangent spaces $T_{x_i} S_{\beta}$ converge to  $\t$, 
then the  Whitney conditions say:
\item{(A)} $T_y S_{\alpha} \subset \t$ 
\item{(B)} $ l \subset \t$.
\end{defin}

\begin{theo}[Stratification, {\cite [p. 1] {G-MP}}]
If $V\subset X$ is a subvariety of a smooth algebraic variety $X$,
then there exists a Whitney stratification for $V$.
\end{theo}

\begin{theo}[Transversality, {\cite [p. 1]{G-MP}}]
Let $Y$ and $X$ be smooth varieties, let $W \subset Y$ and $V \subset X$ 
be their closed subsets with  Whitney stratifications. 
Then the set of maps $f : Y \ra X$ which are transversal to $V$
when restricted to do $W$ is open and dense in the  $\C^\infty$-Whitney topology
on the space of smooth maps from  $Y$ to $X$.
\end{theo}

\bibliography{fakeWeightedProjectiveSpaces}

\bibliographystyle{alpha}

\end{document}